# Nature-Inspired Mateheuristic Algorithms: Success and New Challenges

*Xin-She Yang*
Mathematics and Scientific Computing, National Physical Laboratory, Teddington, TW11 0LW, UK
Present Address: School of Science and Technology, Middlesex University, London NW4 4BT, UK

Many business activities require planning and optimization, this is also true for engineering design, Internet routing, transport scheduling, objective-oriented task management and many other design activities. In fact, optimization is everywhere, the most important part of optimization is the core algorithms used to find optimal solutions to a given problem, though in many cases such algorithms may not exist at all.

## Search Strategy

Let us start by asking a question: suppose you are told there is a treasure in a vast forest, and the treasure, such as a big diamond or a million dollars, is your reward; you are given a limited time, say, a week, to find it. What is your best strategy to find this treasure as quickly as possible?

One strategy is to search the promising areas in the vast forest yourself, inch by inch, and if the treasure is found, it is yours. As the area is vast, it is obviously impossible to go through every square inch of the area within such a limited time. Another strategy is to hire many people/explorers and you promise them to share information and also share any found treasure. Then, if the number of people is large enough, it may be possible to cover the whole area, but the value of the treasure (say, a million dollars) divided among many (say, ten thousands) people many not worth the effort at all. Alternatively, a more sensible approach is probably to use a small group of people, say, twenty or less than a hundred, to search the promising area, share information while searching, and also update and discuss any possible tactics regularly (each hour or each day). This is in fact a swarm intelligence-based approach. As a way to improve the strategy, you as the organizer of the treasure-hunting team, can review the performance of each member of the group regularly, say every day and fire the lazy ones or the least able explorers, and at the same time replace them by recruiting new team members so as to increase the overall performance. This is in fact the selection of the best or elistim.

During the search process, each agent or member of the team can explore the promising region in a quasi-random way. That is, each agent uses the current information to explore any promising area in a local region; if it turns out that no treasure is within this local region, and then the agent can move onto a new, often adjacent, area to do further search. This is the so-called stochastic search strategy. In addition, if each agent is given a walkie talkie, or a mobile phone, to communicate and update their locations and current information, this forms an organized swarm, which may lead to emergent self-organizing behaviour. Imagine a scenario that the team were told that the treasure is hidden at the highest peak of a hilly region in the forest, then they should move towards and climb up the highest peak as quickly as possible; this is essentially a hill-climbing method. If the team were told that the treasure is potentially hidden in a peak but not know which peak, then the group members have to try each possible peak. If they try peak by peak in a sequential manner; that is a hill-climbing with random restart strategy. If they split the group into many small subgroups, then this becomes a parallel hill-climbing strategy. But in reality, there is no such information about the treasure's location, then the best strategy is still yet to be found.

Despite the fact that the best strategy is yet to be found, or may not exist at all, a set of methods have emerged, especially in the last two decades, that they are often surprisingly efficient in practice in solving difficult optimization problems. These methods are called metaheuristic algorithms, and are often nature-inspired, mimicking some successful characteristics in nature. Consequently, these algorithms are also referred to as nature-inspired metaheuristic algorithms. Good examples are Particle Swarm Optimization (PSO), Cuckoo Search (CS) algorithm, Firefly Algorithm (FA), Bat algorithm (BA), Harmony Search (HS), and Ant Colony Optimization (ACO).

The increasing popularity of metaheuristics and swarm intelligence has attracted a great deal of attention in engineering and industry. One of the reasons for this popularity is that nature-inspired metaheuristics are versatile and efficient, and such seemingly simple algorithms can deal with very complex optimisation problems. Metaheuristic algorithms form an important part of contemporary global optimization algorithms, computational intelligence and soft computing.

**Inspiration from Nature**

Nature-inspired algorithms often use multiple interacting agents. A subset of metaheuristcs are often referred to as Swarm Intelligence (SI) based algorithms, and these SI-based algorithms have been developed by mimicking the so-called swarm intelligence characteristics of biological agents such as birds, fish, humans and others. For example, particle swarm optimization was based on the swarming behaviour of birds and fish [1], while the firefly algorithm was based on the flashing pattern of tropical fireflies [2], and cuckoo search algorithm was inspired by the brood parasitism of some cuckoo species.

Nature has been evolving for several hundred million years, and she has found various ingenious solutions to problem-solving and adaption to ever-changing environments. From Darwinian evolution point of view, survival of the fittest will result in the variations and success of species, which can surrive and optimally adapt to environments, and thus selection is a constant pressure that drives the system to improve and adapt for surrival. Any evolutionary advantages over competitors may increase the possibility of reproduction and success of the individuals and the species over the long run.

We can learn from nature by mimicking the successful characteristics of complex systems in nature. Nature-inspired algorithms are still at a very early stage with a relatively short history, comparing with many traditional, well-established methods; however, nature-inspire algorithms have already shown their great potential, flexibility and efficiency with ever-increasing diverse ranges of applications. For example, firefly algorithm was developed by Xin- She Yang in 2008 to mimic the flashing and attraction behaviour of fireflies [2], which leads to a nonlinear dynamical system for optimization using multiple interacting fireflies. Amazingly, firefly algorithm can have some significant advantages over other metaheursistics such as genetic algorithms (GA) and PSO. Two of such advantages are: automatic subgrouping and ability to deal with multimodal problems. Fireflies can automatically subdivide into subgroups and each group can potentially swarm around a local optimum, and all optima (obviously including the global optimum) can be obtained simultaneously if the number of fireflies is much higher than the number of modes. Thus, firefly algorithm can handle multimodal problems very efficiently due to this subgrouping ability. The other advantage is that firefly algorithm does not use velocity, there is no such issues associated with velocities as those in PSO. Consequently, firefly algorithm is much simpler to implement.

As another example, cuckoo search (CS) was developed by Xin- She Yang and Suash Deb in 2009, based on the brooding behavior of some European cuckoo species. By using cuckoo eggs as solutions to an optimization problem, this algorithm produces excellent convergence and high-quality solutions. Enhanced by Lévy flights, cuckoo search can outperform other algorithms such as PSO, GA and ACO for highly nonlinear, global optimization problems. Both cuckoo search and firefly algorithms have been applied in many areas. A quick Google search, at the time of writing in July 2012, leads to about 144 papers on cuckoo search since 2009 and 225 papers on firefly algorithm and their variants since 2008. They certainly form active research topics in optimization and computational intelligence.

New algorithms inspired by natural phenomena, especially biological systems, appear almost every year [1-4]. Even more studies on the extension and improvements on existing algorithms by introducing new components and new applications [3,5]. The literature on these topics is vast, and interested readers can refer to the book by Yang [3] and the references listed in the book.

**Why Metaheuristics?**

New researchers often ask "why metaheuristics?". Indeed, this is a fundamental question to ask in the first steps of solving a given problem. How do we choose the best algorithm and why?

We are often puzzled and often surprised by the excellent efficiency of contemporary nature-inspired algorithms. Seemingly simple algorithms can work 'magic', even for very tough global optimization problems. Many elaborate and sophisticated conventional algorithms often do not work well, despite the fact that conventional algorithms have been well tested for many years. New SI-based metaheuristics often work much better in practice, even though we may not understand why these algorithms actually work. Empirical observations, vast literature and some preliminary convergence analysis all suggest that metaheuristics do work well. Loosely speaking, the success and popularity of metaheuristics can be attributed to the following three factors: algorithm simplicity, ease for implementation, and solution diversity.

Almost all metaheuristic algorithms look simple, and their fundamental characteristics are often derived, directly and indirectly, from nature. Due to the simplicity nature of metaheuristics, they are relatively easy to implement in any programming language. In fact, most algorithms can be coded in fewer than a hundred lines in most programming languages. Such simple algorithms, once implemented properly, can subsequently deal with quite a wide range of optimization problems without much reprogramming.

A key factor may be the balance between solution diversity and solution speed. Ideally, we wish to find the global best solution with the minimum computing effort. For a simple problem, especially a unimodal convex problem, efficient algorithms do exist. For example, conventional algorithms such as hill-climbing or steepest descent methods can find the best solutions in an efficient way for a unimodal problmem. However, real-world problems are not linear, and they certainly are not unimodal. The multimodality and complexity of the problem of interest may mean that we cannot find the global optimality with 100% certainty, unless in a very few limited classes of problems. To reduce the computing time, we often have to sacrifice the diversity of solutions, and consequently, often leading to a local search, or the search process is trapped in a local optimum. In order to escape the local optima, we have to increase the diversity of new solutions so as to potentially reach the true global optimality. The diversity of the solutions in the search process can be achieved in many ways, though randomization and stochastic intervention are often used in most metaheuristics. Now a natural question is how to ensure the proper degree of diversity in the solutions?

In fact, two major components in metaheustics are local exploitation and global exploration. Local exploitation uses local information obtained in local search, and tries to ensure the maximum convergence, while global exploration tends to explore different feasible regions in the whole search space to ensure the global optimality can be achieved with the maximum likelihood. Obviously, these two components are conflicting, and we have to main a tradeoff or balance. For example, bat algorithm was developed by Xin-She Yang in 2010 by using simple rules based on the frequency tuning and echolocation of microbats. Then, this algorithm turns out to be very efficient for a diverse range of problems, and its binary version have been successfully applied to image processing and classifications. The autozooming ability in microbats is manifested in the bat algorithm as automatic adjustment from exploration to exploitation when the global optimality is approaching. This is the first algorithm of its kind in terms of balancing these two key components.

Now the question what is the optimal balance between exploration and exploitation for a given tough optimization task. This is an important question, still without satisfactory answer at the moment. More studies in this area are highly needed.

**Smart Algorithms or Exotic Approaches?**

The popularity of metaheuristics often prompts readers to ask "Can algorithms be intelligent?" The short answer is "possibly" or "it depends".

Artificial intelligence has been an active research area for more than half a century, and new areas such as computational intelligence are going strong as ever. However, unless a Turing test can be really passed in the future, truly intelligent algorithms may be still a long way to go. Obviously, we can define the intelligence by different degrees of mimicking the human intelligence. In that sense, we have been trying to incorporate `intelligence' in the smart algorithms of metaheuristics gradually and incrementally, with some promising results [5].

First, use of memory in the form of selection of the best solutions, elitism and Tabu search is a hint of some intelligence. After all, memory is an important part of human intelligence.

Second, connectionism, interactions and share information can also be considered as `intelligence'. Many algorithms such as artificial neural networks use interactions and connectionism to link inputs to outputs in a complex, implicit manner. In many metaheuristics, multiple agents often can share the best solutions found so far so that new search and solutions are guided by such information.

Thirdly, many algorithms use the so-called swarm intelligence by use certain rules derived from swarm behaviour. These rules essentially ensure the interactions between multiple agents are guided by local information such as the flashing light used in the firefly algorithm or the individual best solution in history found by individual particles in the particle swarm optimization. Mathematically speaking, these interacting agents form biased interacting Markov chains whose convergence rate can be influenced by the structure of the algorithms.

Finally, an algorithm can be called `smart' if it somehow can automatically adjust its behaviour according to the landscape of the objective functions and the information obtained during the search process. If an algorithm with automatic parameter tuning can adjust its algorithm-dependent parameters automatically so as to increase the rate of convergence and reduce the computing cost [5], it may implicitly act in an `intelligent' way.

Obviously, truly intelligent algorithms may only emerge in the far future, however, whatever the forms they may take, they will have a profound impact in almost every area of science, engineering and industrial applications.

**New Challenges**

Despite the increasing popularity of metaheuristics, many crucially important questions remain unanswered. There are two important issues: theoretical framework and the gap between theory and applications. At the moment, the practice of metaheuristics is like heuristic itself, to some extent, by `trial and error'. Mathematical analysis lags far behind, apart from a few, limited, studies on convergence analysis and stability, there is no theoretical framework for analyzing metaheuristic algorithms. I believe mathematical and statistical methods using Markov chains and dynamical systems can be very useful in the future work. There is no doubt that any theoretical progress will provide potentially huge insightful into meteheuristic algorithms.

As there lacks a good theoretical framework, there is thus also a huge gap between theory and applications. Though theory lags behind, applications in contrast are very diverse and active with thousands of papers appearing each year. Accompany this problem, there is another important issue; that is, large-scale problems are yet to be tackled at all. At present, most applications have been tested against toy problems or small-scale benchmarks with a few design variables or at most for problems with a few hundred variables. In reality, many design problems in engineering, business and industry may involve thousands or even millions of variables, we have not seen case studies for such large-scale problems in the literature. In fact, there is no indication that the methods that work for toy benchmarks will work equally well for large-scale problems. In addition to the difference in problem size, there may be some fundamental differences for large-scale problems, and thus the methodology could be significantly different. This still remains a very challenging problem both in theory and in practice.

These important, unresolved, issues also provide golden opportunities for researchers to rethink existing methodology and approaches differently and more fundamentally, and some significant progress may be made in the next ten years. Obviously, any important progress in theory and/or in large-scale pratice will ultimately alter the research landscape in nature-inspired metaheuristics. Maybe, some day, some truly intelligent, self-evolving algorithms may appear to solve tough optimization problems really efficiently.